\input amstex
\voffset 5mm
\magnification=\magstep 1
\documentstyle{amsppt}
\pagewidth{16.6 true cm}
\pageheight{24 true cm}
\TagsOnRight
\emergencystretch=5pt
\tolerance=1000
\def\cl{\mathop{\fam0 cl}}
\def\ri{\mathop{\fam0 ri}}
\def\conv{\mathop{\fam0 conv}}
\def\dom{\mathop{\fam0 dom}}
\def\gr{\mathop{\fam0 gr}}
\def\CL{\mathop{\fam0 CL}}
\topmatter
\title Martingale selection theorem for a stochastic sequence with
       relatively open convex values
\endtitle
\rightheadtext{Martingale selection theorem}
\author    Dmitry B. Rokhlin  \endauthor
\address
\it Dmitry B. Rokhlin\newline
\it \indent Rostov State University, 5 Zorge Street,
Rostov-on-Don~344090, Russia\endaddress
\email\nofrills rokhlin\@math.rsu.ru \endemail
\abstract
For a set-valued stochastic sequence $(G_n)_{n=0}^N$
with relatively open convex values $G_n(\omega)$ we give
a criterion for the existence of an adapted sequence
$(x_n)_{n=0}^N$ of selectors, admitting an equivalent
martingale measure. Mentioned criterion is expressed in terms
of supports of the regular conditional upper distributions of the
elements $G_n$. This result is a refinement of the main result
of author's previous paper (Teor. Veroyatnost. i Primen., 2005, 50:3,
480--500), where the sets $G_n(\omega)$ were assumed to be open
and where were asked if the openness condition can be relaxed.
\endabstract
\endtopmatter
\document
\head Introduction
\endhead
Let $(\Omega,\Cal F,\bold P)$ be a probability space
endowed with the filtration $(\Cal F_n)_{n=0}^N$.
Consider a sequence of $\Cal F_n$-measurable set-valued maps
$\Omega\mapsto G_n(\omega)\subset\Bbb R^d$,
$n=0,\dots,N$ with the nonempty relatively open convex values
$G_n(\omega)$. In this paper we give a criterion
for the existence of a pair, consisting of an adapted single-valued
stochastic process $x=(x_n)_{n=1}^N$,
$x_n(\omega)\in G_n(\omega)$ and a probability measure $\bold Q$
equivalent to $\bold P$ such that $x$ is a martingale under $\bold Q$.
Following [1], we say that {\it the martingale selection problem} (m.s.p.)
is solvable if such a pair $(x,\bold Q)$ exists.

This problem is motivated by some questions of arbitrage theory.
In particular, if the mappings $G_n$ are single-valued, then we
obtain the well-known problem concerning the existence of
an equivalent martingale measure for a given stochastic process $G_n=x_n$.
In this case the solvability of the m.s.p. is equivalent to
the absence of arbitrage in the market, where the discounted
asset price process is described by $x$ [2--5]. It is shown in
[4] that an equivalent martingale measure for $x$ exists iff
the convex hulls of the supports of $x_n-x_{n-1}$ regular conditional
distributions with respect to $\Cal F_{n-1}$ contain the origin as a
point of relative interior [4, Theorem~3,
condition~(g)]. The aim of the present paper is to refine this result.

In the framework of market models with transaction costs [6--8] the
role of equivalent martingale measures is played by
strictly consistent price processes.
This name is assigned to $\bold P$-martingales a.s. taking values in
the relative interior of the random cones
$K^*$, conjugate to the solvency cones $K$. Using the
invariance of $K$ under multiplication, it is easy to show
(see [1, Introduction]) that the existence of a
strictly consistent price process is equivalent to the solvability
of the m.s.p. for the sequence
$(\ri K_n^*)_{n=0}^N$ of relative interiors of $K_n^*$.

In the paper [1] there was obtained a criterion of the solvability of the
m.s.p. under the assumption that the sets $G_n(\omega)$ are open.
This result is not completely satisfactory since, for instance,
it does not include the case of single-valued $G_n$ and it
does not allow the cones $K_n^*$ to have the empty interior.
The last limitation means that the "efficient friction" condition
must be satisfied (according to the terminology of [6]).

In the present paper we refine the main result of [1]
(see Theorem 1). Moreover, the proof given
below, as compared to [1], is considerably simplified.

\head 2. Preliminaries \endhead
Consider a probability space $(\Omega,\Cal F,\bold P)$
and a $\sigma$-algebra $\Cal H\subset\Cal F$.
In the sequel we assume that all $\sigma$-algebras are
complete with respect to $\bold P$ (i.e. they contain all the
subsets of their $\bold P$-negligible sets).
Denote by $\cl A$, $\ri A $, $\conv A$ the closure, the relative
interior, and the convex hull of a subset $A$ of a
finite-dimensional space. Let $\Cal B=\Cal B(\Bbb R^d)$ be the
Borel $\sigma$-algebra of $\Bbb R^d$.

A set-valued map $F$, assigning some set $F(\omega)\subset\Bbb R^d$
to each $\omega\in\Omega$, is called $\Cal H$-measurable if
$\{\omega:F(\omega)\cap V\neq\emptyset\}\in\Cal H$ for any
open set $V\subset\Bbb R^d$.
The graph and the domain of $F$ are defined by
$$ \gr F=\{(\omega,x):x\in F(\omega)\},\ \ \
   \dom F=\{\omega:F(\omega)\neq\emptyset\}.$$
If $\gr F\in\Cal H\otimes\Cal B$, then the mapping
$F$ is $\Cal H$-measurable [9, Corollary II.1.34].

The function $f:\Omega\mapsto\Bbb R^d$ is called a {\it selector}
of a set-valued map $F$ if $f(\omega)\in F(\omega)$
for all $\omega\in\dom F$. Denote by $\Cal S(F,\Cal H)$ the set
of $\Cal H$-measurable selectors of $F$. Note,
that if the set-valued map $F$ is $\Cal H$-measurable, then the
mapping
$$F_*=F I_{\dom F}+\Bbb R^d I_{\Omega\backslash\dom F} \tag 1$$
is also $\Cal H$-measurable and $\Cal S(F,\Cal H)=\Cal S(F_*,\Cal H)$.
Here $I_A(\omega)=1$, $\omega\in A$; $I_A(\omega)=0$,
$\omega\not\in A$.

The countable family $\{f_i\}_{i=1}^\infty$
of an $\Cal H$-measurable selectors is called (an $\Cal H$-measurable)
{\it Castaing representation} for $F$, if the sets
$\{f_i(\omega)\}_{i=1}^{\infty}$ are dense in $F(\omega)$ for all
$\omega\in\dom F$. The set-valued map $F$ with nonempty closed
values is $\Cal H$-measurable iff it admits an $\Cal H$-measurable
Castaing representation [9, Proposition II.2.3].

An element $f\in\Cal S(\conv F,\Cal H)$ is said to have
an $\Cal H$-measurable {\it Caratheodory representation},
if there are some elements
$g_k\in S(F,\Cal H)$, $k=1,\dots,d+1$ and $\Cal H$-measurable
functions
$$\alpha_k\ge 0,\ k=1,\dots,d+1;\ \ \sum_{k=1}^{d+1}\alpha_k=1
\tag 2$$
such that
$f=\sum_{k=1}^{d+1}\alpha_k g_k$ a.s. Under the assumption
$\gr F\in\Cal H\otimes\Cal B$ any element
$f\in\Cal S(\conv F,\Cal H)$ has an $\Cal H$-measurable
Caratheodory representation [10, Theorem 8.2(iii)].

Denote by $\CL=\CL(\Bbb R^d)$ the family of nonempty closed subsets
of $\Bbb R^d$ and let $\Cal E({\CL})$ be the {\it Effros
$\sigma$-algebra}, generated by the sets of the form
$$ A_V=\{D\in{\CL}:D\cap V\neq\emptyset\},$$
where $V$ is an open subset of $\Bbb R^d$.

Suppose $F$ is an $\Cal F$-measurable set-valued map with nonempty
closed values. It follows directly from the definitions that the
corresponding single-valued map
$F:(\Omega,\Cal F) \mapsto (\CL,\Cal E(\CL))$ is measurable.
The measurable space $(\CL,\Cal E(\CL))$ is
a Borel space ([11, Theorem 3.3.10]).
Consequently, the map $F$, considered as a random element taking
values in $(\CL,\Cal E(\CL))$, has the regular conditional
distribution with respect to $\Cal H$
[12, Chapter II, \S 7, Theorem 5].

Thus, there exists a function $\bold P^*:\Omega\times\Cal E(\CL)
\mapsto [0,1]$ with the following properties:
\roster
\item "(i)" for every $\omega$ the function $C\mapsto\bold P^*(\omega,C)$
is a probability measure on $\Cal E(\CL)$;
\item "(ii)" for every $C\in\Cal E(\CL)$ the function
$\omega\mapsto\bold P^*(\omega,C)$ a.s. coincides with the conditional
probability $\bold P(\{F\in C\}|\Cal H)(\omega)$.
\endroster

Following [1], we define the {\it regular conditional upper distribution}
of the mapping $F$ with respect to $\Cal H$ by the formula
$\mu_{F,\Cal H}(\omega,V)=\bold P^*(\omega, A_V)$
for any open subset $V\subset\Bbb R^d$.
The set
$$ \Cal K(F,\Cal H;\omega)=\bigl\{y\in\Bbb R^d:
   \mu_{F,\Cal H}(\omega, \{y':|y'-y|<\varepsilon\})>0\
   \text{for all}\ \varepsilon>0\bigr\}$$
is called the {\it support} of $\mu_{F,\Cal H}(\omega,\cdot)$ [1].
Note, that if $F$ is a single-valued map, then
$\mu_{F,\Cal H}$ is its regular conditional distribution
with respect to $\Cal H$ and $\Cal K(F,\Cal H)$ is the
support of the measure $\mu_{F,\Cal H}$.

The set-valued map $\omega\mapsto\Cal K(F,\Cal H;\omega)$
has nonempty closed values and is $\Cal H$-measurable
[1, Proposition~4(a)]. Let $\{f_i\}_{i=1}^\infty$ be an
$\Cal F$-measurable Castaing representation for $F$. Then
the following equality holds true (see [1, Lemma 1]):
$$ \Cal K(F,\Cal H)=\cl\left(\bigcup_{i=1}^\infty \Cal K(f_i,\Cal H)
\right)\ \text{a.s.} \tag 3$$

If the values of $F$ are empty on a $\bold P$-null set, then
we put $\Cal K(F,\Cal H)=\Cal K(F_*,\Cal H)$,
where $F_*$ is defined by (1). Evidently, equaility (3) still
holds true in this case.

Provided $F(\omega)=\emptyset$ on a set of positive measure,
we put $\Cal K(F,\Cal H)=\emptyset$ for all $\omega$.

\head 3. Main result \endhead
Suppose $\Omega\mapsto G_n(\omega)\subset\Bbb R^d$,
$n=0,1,\dots,N$ is a sequence of $\Cal F_n$-measurable
set-valued maps with nonempty relatively open convex values
$G_n(\omega)$. Define the sequence $(W_n)_{n=0}^N$ of set-valued
maps recursively by
$$ W_N=\cl G_N,$$
$$ W_{n-1}=\cl(G_{n-1}\cap\ri Y_{n-1}),\ \ \
   Y_{n-1}=\conv\Cal K(W_n,\Cal F_{n-1}),\ \ \ 1\le n\le N.$$

This sequence is well-defined and is adapted to the filtration.
Indeed, suppose the map $W_n$ is $\Cal F_n$-measurable.
If $W_n\neq\emptyset$ a.s., then the map
$\conv \Cal K(W_n,\Cal F_{n-1})$ is
$\Cal F_{n-1}$-measurable (see [1, Proposition 4(a)] and
[9, Proposition II.2.26]).
Furthermore, the graphs of the maps $G_{n-1}$, $\ri Y_{n-1}$
are measurable with respect to the $\sigma$-algebra
$\Cal F_{n-1}\otimes\Cal B$ [13, Lemma 1(c)].
Consequently, the map $G_{n-1}\cap\ri Y_{n-1}$
is $\Cal F_{n-1}$-measurable [9, Corollary II.1.34].
Its closure $W_{n-1}$ has the same property [9,
Proposition II.1.8].

Provided $W_n=\emptyset$ on a set of positive measure,
we have $W_{n-1}=\emptyset$ by the definition.
\proclaim{Theorem 1}
The following conditions are equivalent:
\roster
\item "(a)" there exist an adapted to the filtration $(\Cal F_n)_{n=0}^N$
stochastic process $x=(x_n)_{n=0}^N$ and an equivalent to $\bold P$
probability measure $\bold Q$ such that $x_n\in\Cal S(G_n,\Cal F_n)$,
$n\ge 0$ and $x$ is a $\bold Q$-martingale;
\item "(b)" $W_n\neq\emptyset$ a.s., $n=0,\dots,N-1$.
\endroster
\endproclaim

Denote by $\bold E(f|\Cal H)$ the generalized
conditional expectation of the $\Cal F$-measurable
random variable $f$ with respect to $\Cal H$ (under the measure $\bold P$)
[12, p.229], [5, p.117].
The proof of Theorem 1 is based on the following result.

\proclaim{Lemma 1}
Let $F$ be an $\Cal F$-measurable set-valued mapping
with nonempty closed convex values. For any
$\Cal H$-measurable selector $\xi$ of the map
$\ri(\conv \Cal K(F,\Cal H))$ there exist an element
$\eta\in \Cal S(\ri F,\Cal F)$ and an $\Cal F$-measurable
random variable $\gamma>0$ such that
$$ \xi=\bold E(\gamma\eta|\Cal H),\ \ \bold E(\gamma|\Cal H)=1\
\text{a.s.} \tag 4$$
\endproclaim
\demo{Proof}
Let $\{f_i\}_{i=1}^\infty$, $f_i\in\Cal
S(\ri F,\Cal F)$ be a Castaing representation for $\ri F$.
Since $\ri F\in\Cal F\otimes\Cal B$
([13, Lemma 1(c)]), such a representation exists
(see [9, Proposition II.2.17]).

Evidently, $\{f_i\}_{i=1}^\infty$ is also a Castaing representation for
$F$. Applying (3) we get
$$ \xi\in\ri(\conv\Cal K(F,\Cal H))=
\ri\left(\conv\left(\cl\left(\bigcup_{i=1}^\infty
\Cal K(f_i,\Cal H)\right)\right)\right)\ \text{a.s.}$$

Note that for any collection of sets $\{A_i\}_{i=1}^\infty$,
$A_i\subset\Bbb R^d$ the following inclusion holds true
$$ B_1=\ri\left(\conv\left(\cl\left(\bigcup_{i=1}^\infty A_i
\right)\right)\right)\subset\conv\left(\bigcup_{i=1}^\infty
\ri\left(\conv A_i\right)\right)=B_2.$$

Indeed, suppose $x\not\in B_2$. Then by the separation theorem
there exist $p\in\Bbb R^d$, $j\in\Bbb N$,
$\overline y\in A_j$ such that
$$ \langle p,x\rangle\ge\langle p,y\rangle,\ \ \
   y\in A_i,\ i\in\Bbb N;$$
$$ \langle p,x\rangle>\langle p,\overline y\rangle.$$
Here $\langle\cdot,\cdot\rangle$ is the usual scalar
product in $\Bbb R^d$.
Obviously, $\langle p,x\rangle\ge \langle p,z\rangle$ for all
$z\in \cl B_1$. Since $\overline y\in \cl B_1$, it follows that
$\{x\}$ and $\cl B_1$ are properly separated. Therefore,
$x\not\in\ri(\cl B_1)=B_1$.

Putting $A_i=\Cal K(f_i,\Cal H)$, we conclude that
$$ \xi\in\conv\left(\bigcup_{i=1}^\infty\ri
(\conv\Cal K(f_i,\Cal H))\right)\ \text{a.s.}$$
The results of the theory of measurable set-valued maps mentioned
above, readily imply that $\xi$ has an $\Cal H$-measurable
Caratheodory representation:
$$ \xi=\sum_{k=1}^{d+1}\alpha_k\xi_k\ \text{a.s.},\ \
   \xi_k\in\Cal S\left(
   \bigcup_{i=1}^\infty\ri(\conv\Cal K(f_i,\Cal H)),
   \Cal H\right),$$
where $\Cal H$-measurable functions $\alpha_k$ satisfy conditions (2).

Put $A_k^i=\{\omega:\xi_k\in\ri(\conv\Cal K(f_i,\Cal H))\}$
and consider the covering of $\Omega$, consisting of the sets
$A_1^{i_1}\cap\dots\cap A_{d+1}^{i_{d+1}}$, where the upper indexes
run through all natural numbers.
It is easy to show (see [1, Lemma 2]) that there exists an
$\Cal H$-measurable partition $\{D_j\}_{j\in J}$, $J\subset\Bbb N$
of $\Omega$ such that
$$ \emptyset\neq D_j\subset A_1^{i_1}\cap\dots\cap A_{d+1}^{i_{d+1}},
\ \ j\in J,$$
where the set $(i_1,\dots,i_{d+1})$ depends on $j$.

For almost all $\omega\in D_j$ we have
$$ \xi_k\in\ri(\conv \Cal K(f_{{i_k}(j)},\Cal H)),\ \ \
k=1,\dots,d+1,$$
or, in other words,
$ 0\in\ri(\conv\Cal K(\zeta_{kj},\Cal H))\ \text{a.s.},$
where $\zeta_{kj}=I_{D_j}(f_{{i_k}(j)}-\xi_k)$.

According to [4, Theorem 3] it follows that
for any $k\in\{1,\dots,d+1\}$ and $j\in J$
there exists an equivalent to $\bold P$ probability measure $\bold Q_{kj}$
with a.s. bounded density $0<\gamma_{kj}=d\bold Q_{kj}/d\bold P$
such that
$$ I_{D_j}\xi_k=I_{D_j}\bold E_{\bold Q_{kj}}
   (f_{{i_k}(j)}|\Cal H)=\frac{I_{D_j}}{\bold E(\gamma_{kj}|\Cal H)}
   \bold E(\gamma_{kj}f_{{i_k}(j)}|\Cal H)\ \text{a.s.}$$
In the last equality the generalized Bayes formula [5, Ch. V, \S 3a]
is used.

Hence, we get the representation
$$ \xi=\sum_{k=1}^{d+1}\alpha_k\xi_k=\bold E\left
   (\sum_{k=1}^{d+1}\frac{\alpha_k \gamma_{kj}}{\bold E(\gamma_{kj}|\Cal H)}
   f_{{i_k}(j)}\biggr|\Cal H\right)\ \text{a.s.\ on}\ D_j.$$
Here we take into account that the equality
$\bold E(gh|\Cal H)=h\bold E(g|\Cal H)$
holds true if the function $g$ is $\Cal F$-measurable and
$\bold P$-integrable, and the function $h$ is $\Cal H$-measurable
(see the remark in [12, p.~236]).

Put $\beta_{kj}=\gamma_{kj}/\bold E(\gamma_{kj}|\Cal H)$ and
introduce the functions
$$ \gamma_j=\sum_{k=1}^{d+1}\alpha_k\beta_{kj},\ \ \
   \eta_j=\sum_{k=1}^{d+1}\frac{\alpha_k \beta_{kj}}{\gamma_j}
   f_{{i_k}(j)}.$$
We have
$$ \xi=\bold E(\gamma_j\eta_j|\Cal H)\ \text{a.s.\ on\ } D_j.$$
It remains to note that $\gamma_j>0$,
$\bold E(\gamma_j|\Cal H)=1$,
$$ \eta_j\in\conv\{f_{{i_1}(j)},\dots,f_{{i_{d+1}}(j)}\}\subset\ri F\
   \text{a.s.\ on}\ D_j,$$
and the functions
$$\gamma=\sum_{j\in J} I_{D_j} \gamma_j,\ \
\eta=\sum_{j\in J} I_{D_j}\eta_j$$
satisfy conditions (4). The proof of Lemma 1 is complete.
\enddemo
\demo{Proof of Theorem 1} Assume that condition (b) is satisfied.
Starting from an arbitrary selector $x_0\in\Cal
S(\ri W_0,\Cal F_0)$ let us construct adapted sequences
$x_n\in\ri W_n$, $\gamma_n>0$, meeting the conditions
$$ x_{n-1}=\bold E(\gamma_n x_n|\Cal F_n),\ \
   \bold E(\gamma_n|\Cal F_{n-1})=1\ \text{a.s.},\ \ n=1,\dots,N.$$
The existence of the selector $x_0$ is implied by already mentioned
results [13, Lemma 1(c)], [9, Proposition II.2.17].
The existence of the above sequences follows from
Lemma 1, since $x_{n-1}\in\Cal
S(\ri W_{n-1},\Cal F_{n-1})$ imply that $x_{n-1}\in
\Cal S(\ri(\conv \Cal K(W_n,\Cal F_{n-1})),\Cal F_{n-1})$.

Consider the positive $\bold P$-martingale
$$(z_n)_{n=0}^N,\ \ z_0=1,\ \ z_n=\prod_{k=1}^n\gamma_k,\ n\ge 1$$
and the equivalent to $\bold P$ probability measure $\bold Q'$ with
the density $d\bold Q'/d\bold P=z_N$. According to the
generalized Bayes formula we have
$$ x_{n-1}=\frac{1}{z_{n-1}}\bold E(x_n z_n|\Cal F_{n-1})=
   \bold E_{\bold Q'}(x_n|\Cal F_{n-1})\ \text{a.s.}$$
Thus, the process $x$ is a generalized
(or, equivalently, a local) $\bold Q'$-martingale
and it admits an equivalent martingale measure
$\bold Q$ ([4, Theorem 3]).

As long as, moreover, $x_n\in\Cal S(\ri W_n,\Cal F_n)\subset
\Cal S(G_n,\Cal F_n)$, condition (a) is verified.

Now assume that condition (a) is satisfied. Note that $x_N\in G_N\subset
W_N$. Suppose the relations $x_j\in W_j$, $j\ge n$ are already established.
Since
$$0\in\ri(\conv\Cal K(x_n-x_{n-1},\Cal F_{n-1}))\ \text{a.s.},\ n\ge 1$$
(see [4, Theorem 3]) and $\Cal K(x_n,\Cal F_{n-1})\subset
\Cal K(W_n,\Cal F_{n-1})$, it follows that
$$ x_{n-1}\in G_{n-1}\cap\ri(\conv\Cal K(x_n,\Cal F_{n-1}))
\subset G_{n-1}\cap\ri Y_{n-1}\subset W_{n-1}\ \text{a.s.}$$
Particulary, $W_n\neq\emptyset$ a.s. for all $n$. The proof is complete.
\enddemo

In the paper [1] Theorem 1 was proved under one of the following
additional assumptions: (i) the sets $G_n(\omega)$ are open;
(ii) the set $\Omega$ is finite.

\enddocument